\documentclass{IEEEtran}
\usepackage{cite}
\usepackage{amsmath,amssymb,amsfonts}
\usepackage{algorithmic}
\usepackage{graphicx}
\usepackage{textcomp}

\usepackage{subcaption}

\usepackage{amsthm}
\usepackage{mathrsfs}
\usepackage{empheq}
\usepackage{tikz}
\usepackage[utf8]{inputenc}
\usepackage{pgfplots} 
\usepackage{pgfgantt}
\usepackage{pdflscape}
\pgfplotsset{compat=newest} 
\pgfplotsset{plot coordinates/math parser=false}
\usepackage{tikzscale}
\usepackage{standalone}
\usetikzlibrary{shapes,arrows}
\usetikzlibrary{arrows,calc,positioning,fit}

\usepackage{stfloats}

\usepackage[colorinlistoftodos]{todonotes}
\usepackage{mathtools}
\usepackage{enumitem}
\usepackage[T1]{fontenc}

\usepackage{siunitx}
\usepackage{arydshln}

\DeclareMathOperator{\ima}{im}

\newtheorem{theorem}{Theorem}

\newtheorem{assumption}{Assumption}

\newtheorem{lemma}{Lemma}
\newtheorem{corollary}{Corollary}

\newtheorem{example}{Example}

\def\arraystretchval{1.25}

\providecommand{\M}[1]{\begin{bmatrix}#1\end{bmatrix}}
\newcommand{\magne}[1]{#1} 

\def\BibTeX{{\rm B\kern-.05em{\sc i\kern-.025em b}\kern-.08em
    T\kern-.1667em\lower.7ex\hbox{E}\kern-.125emX}}

\begin{document}
\title{Continuity Conditions for \\ Piecewise Quadratic Functions on \\ Simplicial Conic Partitions are Equivalent}
\author{M.~J.~Erlandsen, T.~J.~Meijer, W.~P.~M.~H.~Heemels and S.~J.~A.~M.~van~den~Eijnden
\thanks{The research leading to these results has received funding from the European Research Council under the Advanced ERC Grant Agreement PROACTHIS, no. 101055384.}
\thanks{All authors are with the Control Systems Technology section, Dept. Mechanical Engineering, Eindhoven University of Technology, The Netherlands. E-mail: $\{$m.j.erlandsen, t.j.meijer, m.heemels, s.j.a.m.v.d.eijnden$\}$@tue.nl.}}

\maketitle


\begin{abstract}
Analysis of continuous-time piecewise linear systems based on piecewise quadratic (PWQ) Lyapunov functions typically requires continuity of these functions. Several conditions for guaranteeing continuity of PWQ functions over state-space partitions can be found in the literature.
In this technical note, we show that these continuity conditions are equivalent over so-called simplicial conic partitions. 
As a consequence, the choice of which condition to impose can be based solely on practical considerations such as specific application or numerical aspects, without introducing additional conservatism in the analysis.
\end{abstract}


\section{Introduction}\label{sec:introduction}
Piecewise linear (PWL) systems represent a particular class of switched systems characterised by a partition of the state-space into regions where the system dynamics can be described by linear models \cite{Leenaerts98, Johansson03}. PWL models have become useful within a wide range of applications, including nonsmooth mechanical systems, electrical circuits \cite{Heemels02}, hybrid control \cite{Deenen21, Sharif22, Zaccarian2010, HeeKun_TAC17a}, model predictive control \cite{Bemporad01}, nonlinear system approximation \cite{Camponogara2015}, dynamic optimisation in operations research and
economics \cite{Schumacher2004}, and neural networks \cite{Samanipour24}. 
A particular relevant class of PWL systems are \emph{conewise linear systems} \cite{Camlibel06, Iervolino15}, described by
\begin{equation}\label{eq:conewise}
\dot{x} = A_i x, \textup{ if } Cx \in \mathcal{S}_i,
\end{equation}
where $x \in \mathbb{R}^n$ is the vector of states, $A_i \in \mathbb{R}^{n \times n}$, $i\in \mathcal{N} \coloneq \left\{1,2,\ldots, N\right\}$, and $C \in \mathbb{R}^{m \times n}$, $m \leq n$, are known system matrices, and $\mathcal{S}_i \subseteq \mathbb{R}^m$ are convex polyhedral cones. The collection of polyhedral cones $\mathcal{S}_i$, $i\in \mathcal{N}$, forms a partition of (a subset of) \magne{$\mathbb{R}^n$}.

Stability of conewise linear systems \magne{(and other PWL systems)} as in~\eqref{eq:conewise} is often assessed using PWQ functions of the form 
\begin{equation}\label{eq:V}
V(x) = V_i(x) = x^\top P_i x, \textup{ when } Cx \in \mathcal{S}_i,
\end{equation}
where $P_i = P_i^\top$, 
see, e.g., \cite{Johansson1998, Iervolino15, Yfoulis2004, Deenen21, Ambrosino2015, Eijnden21, Zaccarian2010, HeeKun_TAC17a}.
One reason for their success is that their specific mathematical structure facilitates their construction to be cast into linear matrix inequalities (LMIs), which can be solved systematically using numerical programs. Typical conditions for functions of the form~\eqref{eq:V} to provide a certificate for stability of the conewise linear system in~\eqref{eq:conewise} are formulated in terms of i) positive definiteness of $V(x)$, i.e., $V(0) = 0$ and $V(x) > 0$ for all $x \in \mathbb{R}^n \setminus \left\{0\right\}$, ii) negative definiteness of a suitable generalised time-derivative of $V(x)$, 
and iii) (local Lipschitz) continuity of $V$ over adjacent cones in the partition, that is, $V_i(x) = V_j(x)$ for all $x\in\mathbb{R}^n$ with $Cx \in \mathcal{S}_i \cap \mathcal{S}_j$, $i,j \in \mathcal{N}$. \magne{Note that for PWQ functions as in~\eqref{eq:V}, local Lipschitz continuity is equivalent to continuity \cite[Ch.~4.1.1]{Scholtes2012}.}

Condition i) and ii) are fairly standard and can be guaranteed by searching for matrices $P_i$ that satisfy typical constraints of the form $P_i \succ 0$ and $A_i^\top P_i + P_i A_i \prec 0$, possibly appended with S-procedure relaxation terms \cite{Yakubovich1977} or formulated as a \emph{cone-copositive} problem \cite{Iervolino15}. The arising LMIs can be effectively handled by numerical solvers \cite{YALMIP, MOSEK}. For guaranteeing continuity of the PWQ function over partitions, as stated in condition iii) above, several methods exist in the literature. These methods are either based on posing explicit equality constraints on the matrix $P_i$ \cite{Deenen21, Zaccarian2010}, or on directly incorporating the continuity condition in the parametrisation of the matrix $P_i$ \cite{Johansson1998, Iervolino15, Ambrosino2015}. Both approaches have advantages and disadvantages. For example, equality constraints can be applied to generic partitions but are difficult to solve numerically. The latter results from the fact that solvers work with finite precision and, therefore, return a solution that typically violates the equality constraints \cite{Oehlerking}. Solutions that deal with this numerical inaccuracy have been proposed in specific scenarios \cite{Zaccarian2010,Zaccarian2011}, but the problem remains unsolved in general. 
On the other hand, using a specific matrix parametrisation removes the need for equality constraints, but the parametrisation may be difficult to construct. Despite these apparent differences, it turns out that for certain partitions based on \emph{simplicial} cones, these approaches are \emph{equivalent}. This equivalence result is not completely surprising -- hints at this fact are found scattered across the literature, but has not been proven explicitly and rigorously before. In this technical note, we provide an overview of the various methods for guaranteeing continuity available in the literature, and show their equivalence explicitly. The value of this result lies in demonstrating that none of these approaches introduces additional conservatism in the analysis. Hence, choosing which approach is most suitable can be solely based on practical arguments such as ease of implementation and numerical aspects. 
In addition to the equivalence result, we present a technical lemma inspired by \emph{the non-strict projection lemma} \cite{Meijer2024a}. Although of independent interest, this lemma will be instrumental in proving the aforementioned equivalence.

The remainder of this technical note is organised as follows. Section~\ref{sec:preliminaries} introduces preliminaries on cones and state-space partitions. Section~\ref{sec:lemmas} presents a key lemma. Section~\ref{sec:conditions} discusses continuity conditions for PWQ functions, and Section~\ref{sec:proof} proves their equivalence. Examples are given in Section~\ref{sec:examples}, and conclusions in Section~\ref{sec:conclusion}.

%
%
%
%

\section{Preliminaries}\label{sec:preliminaries}
To make the discussions in this paper precise, in this section, we introduce some mathematical notation used throughout the paper and review some definitions for cones and partitions.

\subsection{Notation}
The set of nonnegative real numbers is denoted $\mathbb{R}_{\geq 0} \coloneq \{ x \in \mathbb{R} \mid x \geq 0 \}$. 
The set of vectors in $\mathbb{R}^n$ (matrices in $\mathbb{R}^{m\times n}$) whose elements are nonnegative real numbers is denoted $\mathbb{R}^n_{\geq 0}$ ($\mathbb{R}^{m\times n}_{\geq 0}$).
The set of symmetric matrices in $\mathbb{R}^{n\times n}$ is denoted $\mathbb{S}^n \coloneq \{ A \in \mathbb{R}^{n\times n} \mid A = A^\top \}$\magne{, and the subset of symmetric matrices with nonnegative entries is denoted $\mathbb{S}^n_{\geq 0}~\coloneq~\{A \in~\mathbb{R}^{n\times n}_{\geq 0} \mid A = A^\top\}$}. 
A positive (semi) definite matrix is denoted $P \succ 0$ ($P \succeq 0$). Similarly, a negative (semi) definite matrix is denoted $P \prec 0$ ($P \preceq 0$).
The transpose of a matrix inverse $\left(A^{-1}\right)^\top$ is compactly written as $A^{-\top}$. The symbol $\star$ is used to complete a symmetric matrix. 
Given a matrix $A\in \mathbb R^{m\times n}$, its image is denoted $\ima A \coloneq \left\{ Av \mid v \in \mathbb{R}^n \right\}$, its kernel is denoted $\ker A \coloneq \left\{ v \in \mathbb{R}^n \mid Av = 0 \right\}$, its (Moore-Penrose) pseudoinverse is denoted $A^+$, and $A_\perp$ denotes any matrix whose columns form a basis of $\ker A$, and thus, $A A_\perp = 0$. The interior of a set $\mathcal{S}$ is denoted $\textup{int}(\mathcal{S})$. 

\subsection{Cones and partitions}\label{ssec:cones_and_partitions}
Given a set of $K$ vectors $z_k \in \mathbb{R}^m$, $k \in \left\{1,2, \ldots, K\right\}$, its \emph{positive hull} (or \emph{conical hull}) is the set of vectors $z \in \mathbb{R}^m$ such that $z = \sum_{k = 1}^K \lambda_k z_k$, with $\lambda_k \geq 0$. A set $\mathcal{C} \subseteq \mathbb{R}^m$ is a \emph{polyhedral cone}, if it is the positive hull of a finite set of vectors. A \emph{face} of a polyhedral cone $\mathcal{C}$, is any set of the form $\mathcal{F} = \mathcal{C} \cap \{ z \in \mathbb{R}^m \mid c^\top z = c_0 \}$, $c_0 \in \mathbb{R}$, that also satisfies $c^\top z \leq c_0$ for all $z \in \mathcal{C}$. A set $\mathcal{S} \subseteq \mathbb{R}^m$ is a \emph{simplicial cone}, if it is the positive hull of $m$ linearly independent vectors.\footnote{The standard definition of a simplicial cone considers the positive hull of (any number of) linearly independent vectors, see, e.g., \cite[Definition~1.39]{Axler}.}

Polyhedral cones with nonempty interior can always be partitioned into a finite number of simplicial cones \cite[Lemma~1.40]{Axler}. For that reason, without loss of generality, in the remainder of this paper, we consider only simplicial cones. Hence, the dynamics of the conewise linear system in~\eqref{eq:conewise} is considered to be defined over simplical cones, that is, the cones $\mathcal{S}_i$, $i \in \mathcal{N}$ in~\eqref{eq:conewise} are assumed to be simplicial. Given a simplicial cone $\mathcal{S} \subseteq \mathbb{R}^m$, there exists a nonsingular matrix $R \in \mathbb{R}^{m \times m}$, such that 
$\mathcal{S} = \{ R\lambda \mid \lambda \in \mathbb{R}^m_{\geq 0} \}$.
The matrix $R$ is called an \emph{extremal ray matrix} of the simplicial cone $\mathcal{S}$. The fact that $R$ is nonsingular follows from our definition of a simplicial cone. The columns of $R$ define the so-called \emph{extremal rays} of the simplicial cone and are uniquely defined up to a positive multiple. 
The set of extremal rays of a simplicial cone $\mathcal{S}$ is denoted $\mathcal{R}_\mathcal{S}$.

Given a set $\mathcal{Z} \subseteq \mathbb{R}^m$ and a positive integer $N$, a \emph{simplicial conic partition} of $\mathcal{Z}$ is a family $\left\{\mathcal{S}_h\right\}_{h = 1}^N$ of simplicial cones satisfying \magne{$\mathcal{Z} = \bigcup_{h = 1}^N \mathcal{S}_h$}, with $\textup{int}(\mathcal{S}_i) \neq \emptyset$ for all $i \in \mathcal{N}$ and $\textup{int}(\mathcal{S}_i) \cap \textup{int}(\mathcal{S}_j) = \emptyset$ for $i,j \in \mathcal{N}, i \neq j$.
\magne{We define the extremal ray matrices $R_i$ of a given simplicial conic partition $\left\{\mathcal{S}_h\right\}_{h = 1}^N$ as follows. Let $r \in \mathbb{N}$ be the number of \emph{distinct} extremal rays of the simplicial conic partition and let $\bar{R} \coloneq \M{r_1 & r_2 & \dots & r_r} \in \mathbb{R}^{m\times r}$, where $r_j$, $j\in \{1,\dots,r\}$ are the distinct extremal rays of the simplicial conic partition. Note that $r\geq m$, since each simplicial cone $\mathcal{S}_i$ is made up of exactly $m$ extremal rays.} For each simplicial cone \magne{$\mathcal{S}_i$} define a so-called extraction matrix $E_i \in \mathbb{R}^{r \times m}$ having its $j$-th row equal to zero for all $r_j \not \in \mathcal{S}_i$, and the remaining rows equal to the rows of the $m$-dimensional identity matrix. Then, the extremal ray matrix of $\mathcal{S}_i$ is given by $R_i = \bar{R}E_i \in \mathbb{R}^{m\times m}$, \magne{see also \cite{Johansson1998, Iervolino15,Rantzer2000}.}

\magne{Throughout this paper we assume the following property of simplicial conic partitions.
\begin{assumption} \label{ass:partition}
    For any two cones of a simplicial conic partition, their intersection $\mathcal{S}_i \cap \mathcal{S}_j$ is a face of both.
\end{assumption}}
\magne{Note that a \emph{face} can be of any dimension, e.g., a single point (the origin). Assumption~\ref{ass:partition} and similar assumptions are fairly standard in the analysis of PWL systems, but are often not addressed explicitly. See, e.g. \cite{Hafstein2017, Camlibel06} in which this assumption is explicitly addressed.} 
A direct result of Assumption~\ref{ass:partition}, is that the extremal rays of the boundary $\mathcal{S}_i \cap \mathcal{S}_j$, are equal to the extremal rays shared by the two cones, i.e., $\mathcal{R}_{\mathcal{S}_i \cap \mathcal{S}_j} = \mathcal{R}_{\mathcal{S}_i} \cap \mathcal{R}_{\mathcal{S}_j}$.
Let the matrix $Z_{ij}$ be a matrix whose columns are equal to the shared extremal rays, i.e., equal to the elements in $\mathcal{R}_{\mathcal{S}_i \cap \mathcal{S}_j}$.
Let $H_{ij} \coloneq ((Z_{ij}^\top)_\perp)^\top$, such that $Z_{ij} = (H_{ij})_\perp$, and thus, $H_{ij}Z_{ij} = 0$.
By definition, the matrix $Z_{ij}$ is tall and has full column rank, whereas the matrix $H_{ij}$ is wide and has full row rank.
Due to the previous definitions and Assumption~\ref{ass:partition}, the boundary shared by two simplicial cones, $\mathcal{S}_i \cap \mathcal{S}_j$, satisfies
\begin{equation}\label{eq:2_boundary_representation}
    \mathcal{S}_i \cap \mathcal{S}_j =  \{ Z_{ij} v \mid v \geq 0 \} \subseteq \ima Z_{ij} = \ker H_{ij}.
\end{equation} 

\section{Technical Lemma}\label{sec:lemmas}
In this section, we present a technical lemma in the spirit of the non-strict projection lemma in~\cite{Meijer2024a}. Although this lemma is of independent interest, it will be useful in proving equivalence of the continuity conditions in Section~\ref{sec:conditions}.

\begin{lemma}\label{lem:equality-projection}
    Let $\boldsymbol{U}\in\mathbb{R}^{m\times n}$, $\boldsymbol{V}\in\mathbb{R}^{p\times n}$ and $\boldsymbol{Q}\in\mathbb{S}^n$. Consider the following statements:
    \begin{enumerate}[labelindent=5pt,labelwidth=\widthof{\ref{item:P3}},itemindent=0em,leftmargin=!,label=(L\ref{lem:equality-projection}.\arabic*),topsep=6pt]
        \item\label{item:P1} There exists a matrix $\boldsymbol{X}\in\mathbb{R}^{m\times p}$ such that 
        $$\boldsymbol{Q}+\boldsymbol{U}^\top \boldsymbol{XV} + \boldsymbol{V}^\top \boldsymbol{X}^\top \boldsymbol{U}= 0;$$
        \item\label{item:P2}$x^\top \boldsymbol{Q} x = 0 \text{ for all } x \in \ker \boldsymbol{U} \cup \ker \boldsymbol{V};$
        \item\label{item:P3} $\ker \boldsymbol{U}\cap\ker \boldsymbol{V} \subseteq\ker \boldsymbol{Q}$.
    \end{enumerate}
    Then,~\ref{item:P1} holds if and only if~\ref{item:P2} and~\ref{item:P3} hold.
\end{lemma}
\begin{proof}
    {\bf Necessity:} Suppose that~\ref{item:P1} holds.
    Then, using the fact that either $\boldsymbol{U}x=0$ or $\boldsymbol{V}x=0$ when $x \in \ker \boldsymbol{U} \cup \ker \boldsymbol{V}$, it follows that~\ref{item:P2} holds.
    Due to~\ref{item:P1}, it holds, for any $x\in\mathbb{R}^n$, that
    \begin{equation}\label{eq:zero2}
        \left(\boldsymbol{Q}+\boldsymbol{U}^\top \boldsymbol{XV}+\boldsymbol{V}^\top \boldsymbol{X}^\top \boldsymbol{U}\right)x = 0.
    \end{equation}
    Let $x\in\ker \boldsymbol{U}\cap \ker \boldsymbol{V}$. Then, it holds that
    \begin{equation}
        \left(\boldsymbol{Q}+\boldsymbol{U}^\top \boldsymbol{XV}+\boldsymbol{V}^\top \boldsymbol{X}^\top \boldsymbol{U}\right)x = \boldsymbol{Q}x \stackrel{\eqref{eq:zero2}}{=}0,
    \end{equation}
    i.e.,~\ref{item:P3} holds.

    {\bf Sufficiency:} Suppose~\ref{item:P2} and~\ref{item:P3} hold. Let $\boldsymbol{T}\in\mathbb{R}^{n\times n}$ be a nonsingular matrix, whose columns in the partition $\boldsymbol{T}=\begin{bmatrix} \boldsymbol{T}_1 & \boldsymbol{T}_2 & \boldsymbol{T}_3 & \boldsymbol{T}_4\end{bmatrix}$ are chosen to satisfy
    \begin{align}
        \ima \begin{bmatrix}
            \boldsymbol{T}_1 & \boldsymbol{T}_3
        \end{bmatrix} &= \ker \boldsymbol{U},\label{eq:T1T3}\\
        \ima \begin{bmatrix}
            \boldsymbol{T_2} & \boldsymbol{T}_3
        \end{bmatrix} &= \ker \boldsymbol{V},\label{eq:T2T3}\\
        \ima \boldsymbol{T}_3 &= \ker \boldsymbol{U}\cap \ker \boldsymbol{V}.\label{eq:T3}
    \end{align}
    Clearly,~\ref{item:P1} is equivalent to the existence of $\boldsymbol{X}\in\mathbb{R}^{m\times p}$ such that
    \begin{equation}\label{eq:Y}
        \boldsymbol{Y}\coloneqq \boldsymbol{T}^\top\left(\boldsymbol{Q}+\boldsymbol{U}^\top \boldsymbol{XV}+\boldsymbol{V}^\top \boldsymbol{X}^\top \boldsymbol{U}\right)\boldsymbol{T}=0.
    \end{equation}
    We partition $\boldsymbol{W}\coloneqq \boldsymbol{T}^\top \boldsymbol{QT}$ in accordance with $\boldsymbol{T}$ to obtain
    \begin{equation}
        \boldsymbol{W} = \boldsymbol{W}^\top = \begin{bmatrix}
            \boldsymbol{W}_{11} & \boldsymbol{W}_{12} & \boldsymbol{W}_{13} & \boldsymbol{W}_{14}\\
            \star & \boldsymbol{W}_{22} & \boldsymbol{W}_{23} & \boldsymbol{W}_{24}\\
            \star & \star & \boldsymbol{W}_{33} & \boldsymbol{W}_{34}\\
            \star & \star & \star & \boldsymbol{W}_{44}
        \end{bmatrix}.\label{eq:W}
    \end{equation}
    Using~\eqref{eq:T1T3},~\eqref{eq:T2T3} and~\eqref{eq:T3}, we write the term $(\boldsymbol{UT})^\top \boldsymbol{X}(\boldsymbol{VT})$ in~\eqref{eq:Y} as
    \begin{equation}
        \begin{bmatrix}
            \boldsymbol{U}\boldsymbol{T}_2 & \boldsymbol{U}\boldsymbol{T}_4
        \end{bmatrix}^\top \boldsymbol{X} \begin{bmatrix}
            \boldsymbol{V}\boldsymbol{T}_1 & \boldsymbol{V}\boldsymbol{T}_4
        \end{bmatrix}\eqqcolon \begin{bmatrix}
            \boldsymbol{K} & \boldsymbol{L}\\
            \boldsymbol{M} & \boldsymbol{N}
        \end{bmatrix},\label{eq:UXV}
    \end{equation}
    where, due to~\eqref{eq:T1T3} and~\eqref{eq:T2T3}, $\begin{bmatrix} \boldsymbol{U}\boldsymbol{T}_2 & \boldsymbol{U}\boldsymbol{T}_4\end{bmatrix}$ and $\begin{bmatrix} \boldsymbol{V}\boldsymbol{T}_1 & \boldsymbol{V}\boldsymbol{T}_4\end{bmatrix}$ have full column rank. Hence, using~\eqref{eq:W} and~\eqref{eq:UXV},~\eqref{eq:Y} reads as
    \begin{align}
        \boldsymbol{Y}&=\boldsymbol{Y}^\top =\def\arraystretch{\arraystretchval}\left[\begin{array}{@{}c;{2pt/2pt}c@{}}
            \boldsymbol{Y}_1 & \boldsymbol{Y}_2\\\hdashline[2pt/2pt]
            \star & \boldsymbol{Y}_3
        \end{array}\right] \nonumber \\
        &= \def\arraystretch{\arraystretchval}\left[\begin{array}{@{}ccc;{2pt/2pt}c@{}}
            \boldsymbol{W}_{11} & \boldsymbol{W}_{12} + \boldsymbol{K}^\top & \boldsymbol{W}_{13} & \boldsymbol{W}_{14}+\boldsymbol{M}^\top\\
            \star & \boldsymbol{W}_{22} & \boldsymbol{W}_{23} & \boldsymbol{W}_{24}+\boldsymbol{L}\\
            \star & \star & \boldsymbol{W}_{33} & \boldsymbol{W}_{34}\\\hdashline[2pt/2pt]
            \star & \star & \star & \boldsymbol{W}_{44} + \boldsymbol{N}+\boldsymbol{N}^\top
        \end{array}\right] \nonumber  \\
        &= 0.
    \end{align}
    It follows from~\ref{item:P2} that
    \begin{equation}
        \begin{bmatrix}
            \boldsymbol{W}_{11} & \boldsymbol{W}_{13}\\
            \star & \boldsymbol{W}_{33}
        \end{bmatrix} = 0\text{ and }\begin{bmatrix}
            \boldsymbol{W}_{22} & \boldsymbol{W}_{23}\\
            \star & \boldsymbol{W}_{33}
        \end{bmatrix} = 0.
    \end{equation}
    Clearly, to ensure that $\boldsymbol{Y}_1=0$, we should construct $\boldsymbol{X}$ such that $\boldsymbol{K} = -\boldsymbol{W}_{12}^\top$. Similarly, we will aim to construct $\boldsymbol{X}$ such that $\boldsymbol{L}=-\boldsymbol{W}_{24}$, $\boldsymbol{M}=-\boldsymbol{W}_{14}^\top$ and $\boldsymbol{N}=-\frac{1}{2}\boldsymbol{W}_{44}$. Note that, due to $\begin{bmatrix} \boldsymbol{U}\boldsymbol{T}_2 & \boldsymbol{U}\boldsymbol{T}_4\end{bmatrix}$ and $\begin{bmatrix} \boldsymbol{V}\boldsymbol{T}_1 & \boldsymbol{V}\boldsymbol{T}_4\end{bmatrix}$ having full column rank, we can construct such $\boldsymbol{X}$ by taking
    \begin{equation}
    \begin{split}
        \boldsymbol{X} &= \begin{bmatrix}
            \left(\boldsymbol{U}\boldsymbol{T}_2\right)^\top\\
            \left(\boldsymbol{U}\boldsymbol{T}_4\right)^\top
        \end{bmatrix}^+\begin{bmatrix}
            \boldsymbol{K} & \boldsymbol{L}\\
            \boldsymbol{M} & \boldsymbol{N}
        \end{bmatrix}\begin{bmatrix}
            \boldsymbol{V}\boldsymbol{T}_1 & \boldsymbol{V}\boldsymbol{T}_4
        \end{bmatrix}^+ \\
        &= \begin{bmatrix}
            \left(\boldsymbol{U}\boldsymbol{T}_2\right)^\top\\
            \left(\boldsymbol{U}\boldsymbol{T}_4\right)^\top
        \end{bmatrix}^+\begin{bmatrix}
            -\boldsymbol{W}_{12}^\top & -\boldsymbol{W}_{24}\\
            -\boldsymbol{W}_{14}^\top & -\frac{1}{2}\boldsymbol{W}_{44}
        \end{bmatrix}\begin{bmatrix}
            \boldsymbol{V}\boldsymbol{T}_1 & \boldsymbol{V}\boldsymbol{T}_4
        \end{bmatrix}^+.
    \end{split}
    \end{equation}
    Note that all entries now equal zero except for $\boldsymbol{W}_{34}$. Hence, it remains to show that $\boldsymbol{W}_{34}=0$. It follows from~\ref{item:P3} that $\boldsymbol{Q}\boldsymbol{T}_3=0$ and, thus, $\boldsymbol{W}_{34}^\top = \boldsymbol{T}_4^\top \boldsymbol{Q}\boldsymbol{T}_3 = 0$.
\end{proof}

\noindent Lemma~\ref{lem:equality-projection} is closely related to the non-strict projection lemma \cite{Meijer2024a}, but it deals with equalities instead of (matrix) inequalities. Interestingly, it turns out that, as in the non-strict projection lemma, an additional coupling condition~\ref{item:P3} is needed to achieve the equivalence in Lemma~\ref{lem:equality-projection}.

Next, we introduce two useful corollaries of Lemma~\ref{lem:equality-projection}.

\begin{corollary}\label{cor:onesided_projection}
    Let $\boldsymbol{U}\in\mathbb{R}^{m\times n}$ and let $\boldsymbol{Q}\in\mathbb{S}^n$. \magne{The following two statements are equivalent:}
    \begin{enumerate}[labelindent=5pt,labelwidth=\widthof{\ref{item:P3}},itemindent=0em,leftmargin=!,label={(C\ref{cor:onesided_projection}.\arabic*}),topsep=6pt]
        \item\label{item:OP1} There exists a matrix $\boldsymbol{X}\in\mathbb{R}^{m\times n}$ such that 
        $$\boldsymbol{Q}+\boldsymbol{U}^\top \boldsymbol{X} + \boldsymbol{X}^\top \boldsymbol{U}= 0;$$
        \item\label{item:OP2} $x^\top \boldsymbol{Q}x=0$ for all $x \in \ker \boldsymbol{U}$.
    \end{enumerate}
\end{corollary}
Corollary~\ref{cor:onesided_projection} follows from Lemma~\ref{lem:equality-projection} with $\boldsymbol{V}\coloneq I_n$. To see this, note that $\ker \boldsymbol{V} = \{0\}$, and thus,~\ref{item:P3} trivially holds.

\begin{corollary}\label{cor:equality_finsler}
    Let $\boldsymbol{U}\in\mathbb{R}^{m\times n}$ and let $\boldsymbol{Q}\in\mathbb{S}^n$. \magne{The following two statements are equivalent:}
    \begin{enumerate}[labelindent=5pt,labelwidth=\widthof{\ref{item:F1}},itemindent=0em,leftmargin=!,label={(C\ref{cor:equality_finsler}.\arabic*}),topsep=6pt]
        \item\label{item:F1} There exists a symmetric matrix $\boldsymbol{X}\in\mathbb{S}^{m}$ such that 
        $$\boldsymbol{Q}+\boldsymbol{U}^\top \boldsymbol{XU}= 0;$$
        \item\label{item:F2} $\ker \boldsymbol{U} \subseteq\ker \boldsymbol{Q}$.
    \end{enumerate}
\end{corollary}
Corollary~\ref{cor:equality_finsler} follows from Lemma~\ref{lem:equality-projection} with $\boldsymbol{V} \coloneq \frac 1 2 \boldsymbol{U}$, in which case $\ker \boldsymbol{U} = \ker \boldsymbol{V}$. Thus,~\ref{item:P3} simplifies to~\ref{item:F2}, which immediately implies~\ref{item:P2}. Corollary~\ref{cor:equality_finsler} is closely related to the non-strict Finsler's lemma \cite{Meijer2024b}, but it deals with equalities instead of (matrix) inequalities.

\section{Continuity Conditions}\label{sec:conditions}
In this section, we formalise the equivalence of different conditions that can be found in the literature for guaranteeing continuity of a PWQ function. In particular, we consider PWQ functions of the form as in~\eqref{eq:V}, i.e., where $P_i \in \mathbb{S}^n$, $i \in \mathcal{N}$, $C\in \mathbb{R}^{m\times n}$ has full row rank, and $\mathcal{S}_i$ are simplicial cones.
We want to guarantee continuity of these PWQ functions and thus, local Lipschitz continuity, in order for them to be useful in stability analysis.

Before stating the main theorem, let us emphasise that definitions from Section~\ref{sec:preliminaries} are used, e.g., for the matrices $Z_{ij}$, $H_{ij}$, $E_i$, and $R_i$.

\begin{theorem}\label{thm:1}
Let $\mathcal{N} \coloneq \{1,2, \ldots, N\}$. Consider a simplicial conic partition $\left\{\mathcal{S}_i\right\}_{i=1}^N$ of a set $\mathcal{Z} \subseteq \mathbb{R}^m$, and a set $\{P_i\}_{i=1}^N$ of symmetric matrices $P_i \in \mathbb{S}^n$, $i\in\mathcal{N}$. Then, the following statements are equivalent:
\begin{enumerate}[labelindent=5pt,labelwidth=\widthof{\ref{item:C2}},itemindent=0em,leftmargin=!,label={(T\ref{thm:1}.\arabic*}),topsep=6pt]
    \item\label{item:C1} The matrices $P_1, P_2, \ldots, P_N$ satisfy, for all $i,j\in\mathcal{N}$
\begin{equation}\label{eq:C1}
    x^\top (P_i-P_j) x = 0, \textup{ for all } Cx \in \mathcal{S}_i \cap \mathcal{S}_j,
\end{equation}
    and thus, the function $V(x)$ as in~\eqref{eq:V} is continuous.
    \item\label{item:C12_temp} The matrices $P_1, P_2, \ldots, P_N$ satisfy, for all $i,j\in\mathcal{N}$
\begin{equation}\label{eq:C12_temp}
    x^\top (P_i-P_j) x = 0, \textup{ for all } Cx \in \ima Z_{ij}. 
\end{equation}
    \item\label{item:C2} Let 
    \begin{equation*}
    W_{ij} = T^{-1} \M{Z_{ij} & 0 \\ 0 & I}, \textup{ with } T \coloneq \M{C \\ C_\perp^\top} .
\end{equation*}
    For all $i,j\in\mathcal{N}$, it holds that
\begin{equation}\label{eq:C2}
    W_{ij}^\top (P_i - P_j)W_{ij} = 0.
\end{equation}
    \item\label{item:C3} Let
    \begin{equation*}
    F_i = \M{E_i R_i^{-1}C \\ V},
\end{equation*}
    where $V$ is any matrix 
    that satisfies \magne{$\ima V^\top \supseteq \ker C$}.
    There exists a symmetric matrix $\Phi$, such that, for all $i\in\mathcal{N}$
    \begin{equation}\label{eq:C3}
    P_i = F_i^\top \Phi F_i.
\end{equation}
    \item\label{item:C4} There exist matrices $\Gamma_{ij}$, for all $i,j\in\mathcal{N}$, such that
\begin{equation}\label{eq:C4}
    P_i-P_j + (H_{ij}C)^\top \Gamma_{ij} + \Gamma_{ij}^\top (H_{ij}C) = 0.
\end{equation}
\end{enumerate}
\end{theorem}

\begin{figure}
    \centering
    \includegraphics[width=0.9\linewidth]{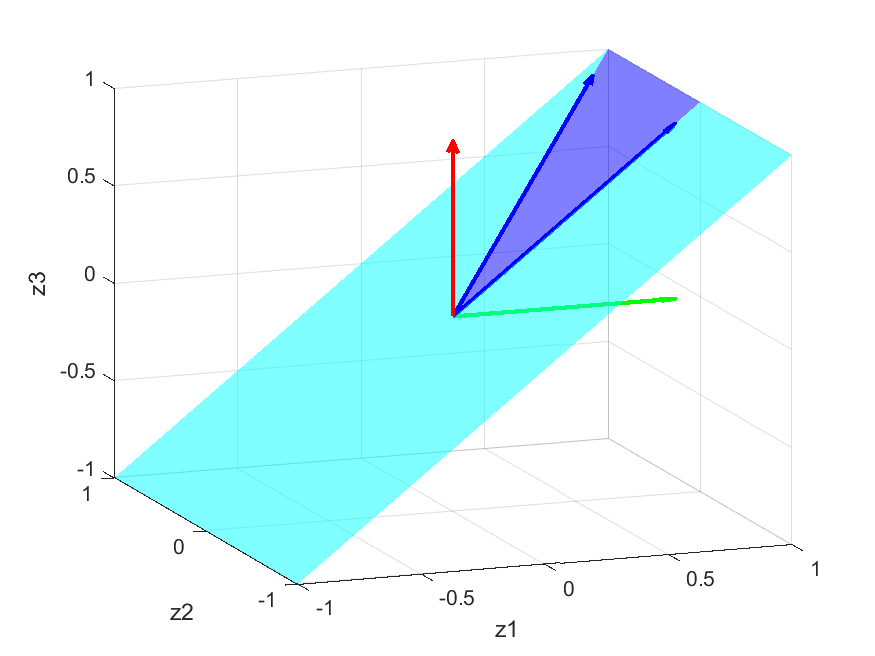}
    \caption{Example illustrating the difference between the regions considered in the first two items of Theorem~\ref{thm:1}. The figure includes the distinct extremal rays of two example cones, $\mathcal{S}_1$ and $\mathcal{S}_2$ (in red and green), their shared extremal rays (in blue), their boundary region, $\mathcal{S}_1 \cap \mathcal{S}_2$ (also in blue), and the extended boundary region, $\ima Z_{ij}$ (in cyan).}
    \label{fig:regions}
\end{figure}

In the next section, we will give an explicit proof of the equivalence in Theorem~\ref{thm:1}. However, before continuing with the proof, we provide a few comments and discussions on the various elements of Theorem~\ref{thm:1}:

\begin{enumerate}
\item Item~\ref{item:C1} expresses necessary and sufficient conditions for continuity of a PWQ function as in \eqref{eq:V}, over generic state-space partitions. However, we will only show its equivalence with the other conditions,~\ref{item:C12_temp}-\ref{item:C4}, over simplicial conic partitions. Hence, in general, equivalence may not be guaranteed.
Conditions~\ref{item:C1} and~\ref{item:C12_temp} require checking an infinite number of equalities, that is, one for each $x \in \mathbb{R}^n$. On the other hand, \ref{item:C2}--\ref{item:C4} express continuity conditions in terms of computationally tractable conditions on the matrices $P_i$ directly. \magne{Note that~\ref{item:C2} expresses the values of $x$ on the boundary in terms of so-called basis functions given by $W_{ij}$.}
\item An example illustrating the difference between the sets $\mathcal{S}_i \cap \mathcal{S}_j$ and $\ima Z_{ij}$, used in~\ref{item:C1} and~\ref{item:C12_temp}, is shown in Fig.~\ref{fig:regions}. The set $\ima Z_{ij}$ is the minimal linear subspace of $\mathbb{R}^m$ that contains $\mathcal{S}_i \cap \mathcal{S}_j$ (minimal in the sense that its dimension is equal to the dimension of $\mathcal{S}_i \cap \mathcal{S}_j$, or equivalently, that $\ima Z_{ij}$ is equal to the intersection of all possible linear subspaces of $\mathbb{R}^m$ that contain $\mathcal{S}_i \cap \mathcal{S}_j$). We can say that $\ima Z_{ij}$ is the so-called \emph{linear hull} of the boundary region $\mathcal{S}_i \cap \mathcal{S}_j$.  Surprisingly, the equivalence between \ref{item:C1} and~\ref{item:C12_temp} means that, for PWQ functions, continuity on the boundary, $\mathcal{S}_i \cap \mathcal{S}_j$, is equivalent to continuity on the whole (generalised) plane (of some dimension) containing $\mathcal{S}_i \cap \mathcal{S}_j$.
\item Continuity conditions of the form presented in \ref{item:C2} are used in, e.g., \cite{Deenen21}, \cite{Zaccarian2010}. The equality constraint in~\eqref{eq:C2} is simple to formulate, but generally difficult to solve numerically. The reason for this, is that solvers work with finite precision and, as a result, often return solutions that slightly violate the equality constraints (see \cite[Section~4.5.2]{Oehlerking}).
\magne{In~\cite{Zaccarian2010}, such violations were dealt with for planar partitions only (i.e., $m=2$). 
If one instead wishes to verify whether a given PWQ function is continuous, the conditions \ref{item:C2} are convenient, as they simply involve checking whether a number of equalities are satisfied.
}
Note that the matrix $T = \M{C^\top & C_\perp}^\top$ is nonsingular, since $C$ is assumed, without loss of generality, to have full row rank.
\magne{When $C$ also has full column rank, $C_\perp$ is omitted, and $W_{ij}=C^{-1}Z_{ij}$.}
\item The parametrisation \magne{$P_i = F_i^\top \Phi F_i$} in~\ref{item:C3} was first proposed in \cite{Johansson1998}, and has been used successfully, e.g. in \cite{Iervolino15, Ambrosino2015}. This parametrisation removes the need for explicit equality constraints, which may provide a significant advantage from a computational point of view \cite{Eijnden21}. For the matrix $V$, a natural choice is $V=C_\perp^\top$. This choice minimises the number of parameters in the matrix $\Phi$, which can be numerically beneficial. Another simple choice is $V=I$, which avoids the need to compute the matrix $C_\perp$, but leads to more parameters in $\Phi$. 
\magne{Note that the requirement $\ima V^\top \supseteq \ker C$ is equivalent to the matrix $\M{C^\top & V^\top}^\top$ having full column rank. Consequently, when the matrix $C$ has full column rank, $V$ is omitted, and $F_i = E_i R_i^{-1}C$. This construction of the so-called continuity matrices $F_i$ coincides with that proposed in~\cite[Section~\uppercase\expandafter{\romannumeral6\relax}]{Rantzer2000}, although the latter considers the more general class of simplex partitions. 
When $C$, on the other hand, is a wide matrix with $m<n$ linearly independent rows, the continuity conditions in \cite{Rantzer2000} are, to the best of our knowledge, no longer necessary for continuity of \eqref{eq:V}, but only sufficient.} 
\item \magne{Although outside the scope of this paper, algorithms for partitioning a set into smaller regions play an important role in reducing conservatism (by enabling more versatile PWQ functions); see, e.g.,~\cite{Iervolino15,Ambrosino2015,Hafstein2017}.}
\end{enumerate}

\section{Proof of Theorem~\ref{thm:1}}\label{sec:proof}	
In this section, we present the proof of Theorem~\ref{thm:1}. The proof is carried out in the following order: $\ref{item:C1} \iff~\ref{item:C12_temp}$, $\ref{item:C12_temp} \iff~\ref{item:C2}$, $\ref{item:C12_temp} \iff~\ref{item:C4}$, $\ref{item:C1} \implies~\ref{item:C3}$, $\ref{item:C3} \implies~\ref{item:C1}$.

$\mathbf{\ref{item:C1} \iff~\ref{item:C12_temp}}$. Since $\mathcal{S}_i \cap \mathcal{S}_j \subseteq \ima Z_{ij}$, the necessity, $\ref{item:C1} \impliedby~\ref{item:C12_temp}$, is trivial. Hence, the focus is on the sufficiency, $\ref{item:C1} \implies~\ref{item:C12_temp}$.
Suppose that $x^\top (P_i-P_j) x = 0$ for all $Cx \in \mathcal{S}_i \cap \mathcal{S}_j$. Consider the state transformation $\bar{x} = T x$, with $T = \M{C^\top & C_\perp}^\top$, such that $\bar{x} = \M{z^\top & \hat{x}^\top}^\top$, with $z = Cx \in \mathbb{R}^m$ and $\hat{x} = C_\perp^\top x \in \mathbb{R}^{n-m}$.
Now, partition the matrix $\bar{P}_i \coloneq T^{-\top } P_i T^{-1}$ as
\begin{equation}\label{eq:Pipart_new}
    \bar{P}_i = \bar{P}_i^\top = \M{\bar{P}_{i}^{11} & \star \\ \bar{P}_{i}^{21} & \bar{P}_{i}^{22}},
\end{equation}
according to $(z, \hat{x})$, such that 
\begin{align}
    0 &= x^\top (P_i-P_j) x = \bar{x}^\top \left( \bar{P}_i - \bar{P}_j \right) \bar{x} \nonumber \\
    &= \M{z \\ \hat{x}}^\top \left( \M{\bar{P}_{i}^{11} & \star \\ \bar{P}_{i}^{21} & \bar{P}_{i}^{22}} - \M{\bar{P}_{j}^{11} & \star \\ \bar{P}_{j}^{21} & \bar{P}_{j}^{22}}  \right) \M{z \\ \hat{x}} \nonumber \\
    &= z^\top (\bar{P}_i^{11}-\bar{P}_j^{11})z + 2 \hat{x}^\top (\bar{P}_i^{21}-\bar{P}_j^{21})z \nonumber \\
    &\quad + \hat{x}^\top (\bar{P}_i^{22}-\bar{P}_j^{22})\hat{x} \label{eq:cont_bar},
\end{align}
for all $z \in \mathcal{S}_i \cap \mathcal{S}_j$ and all $\hat x \in \mathbb{R}^{n-m}$. Since $0 \in \mathcal{S}_i \cap \mathcal{S}_j$,~\eqref{eq:cont_bar} may be evaluated separately for $z=0$ and $\hat{x} = 0$. As such, one finds that
\begin{subequations} \label{eq:Pbari_cont}
\begin{align}
    z^\top (\bar{P}_i^{11}-\bar{P}_j^{11})z &= 0,  \textup{ for all } z\in\mathcal{S}_i \cap \mathcal{S}_j , \label{eq:Pbari11_cont} \\
    \hat{x}^\top (\bar{P}_i^{21}-\bar{P}_j^{21})z &= 0, \textup{ for all } z\in\mathcal{S}_i \cap \mathcal{S}_j , \; \hat{x} \in \mathbb{R}^{n-m},\label{eq:Pbari21_cont}\\
    \hat{x}^\top (\bar{P}_i^{22}-\bar{P}_j^{22})\hat{x} &= 0, \textup{ for all } \hat{x} \in \mathbb{R}^{n-m} . \label{eq:Pbari22_cont}
\end{align}
\end{subequations}
Firstly,~\eqref{eq:Pbari22_cont} implies, due to symmetry of $\bar{P}_i^{22}-\bar{P}_j^{22}$, that
\begin{equation} \label{eq:Pbari22_gen}
    \bar{P}_i^{22}-\bar{P}_j^{22} = 0.
\end{equation}
Secondly, since $\mathcal{R}_{\mathcal{S}_i \cap \mathcal{S}_j} \subseteq \mathcal{S}_i \cap \mathcal{S}_j$, it follows from~\eqref{eq:Pbari21_cont} that $\hat{x}^\top (\bar{P}_i^{21}-\bar{P}_j^{21})r = 0$ for every $r \in \mathcal{R}_{\mathcal{S}_i \cap \mathcal{S}_j}$ and all $\hat{x} \in \mathbb{R}^{n-m}$. As such,
\begin{equation} \label{eq:Pbari21_gen}
    (\bar{P}_i^{21}-\bar{P}_j^{21})r = 0, \textup{ for every } r \in \mathcal{R}_{\mathcal{S}_i \cap \mathcal{S}_j}.
\end{equation}
Thirdly, it follows from~\eqref{eq:Pbari11_cont} that
\begin{equation} \label{eq:Pbari11_gen1}
    r^\top (\bar{P}_i^{11}-\bar{P}_j^{11})r = 0, \textup{ for every } r\in \mathcal{R}_{\mathcal{S}_i \cap \mathcal{S}_j}. 
\end{equation}
Note that
\begin{equation} \label{eq:rm_sum}
    r_m + r_n \subseteq \mathcal{S}_i \cap \mathcal{S}_j, \textup{ for every } r_m,r_n \in \mathcal{R}_{\mathcal{S}_i \cap \mathcal{S}_j},
\end{equation}
that is, the sum of two shared extremal rays is contained in the boundary region $\mathcal{S}_i \cap \mathcal{S}_j$. As such, by substituting $z = r_m + r_n$ into~\eqref{eq:Pbari11_cont}, it follows that, for every $r_m,r_n \in \mathcal{R}_{\mathcal{S}_i \cap \mathcal{S}_j}$
\begin{align}
    0 &\stackrel{\eqref{eq:Pbari11_cont}}{=} z^\top (\bar{P}_i^{11}-\bar{P}_j^{11})z \nonumber \\
    &= \left(r_m+r_n \right)^\top (\bar{P}_{i}^{11} - \bar{P}_{j}^{11}) \left(r_m+r_n\right) \nonumber \\
    &= r_m^\top (\bar{P}_{i}^{11} - \bar{P}_{j}^{11}) r_m + r_n^\top (\bar{P}_{i}^{11} - \bar{P}_{j}^{11}) r_n \nonumber \\
    &\quad+ 2 r_m^\top (\bar{P}_{i}^{11} - \bar{P}_{j}^{11}) r_n \nonumber \\
    &\stackrel{\eqref{eq:Pbari11_gen1}}{=} 2 r_m^\top (\bar{P}_{i}^{11} - \bar{P}_{j}^{11}) r_n . \label{eq:Pbari_11_gen2}
\end{align}
Note that~\eqref{eq:Pbari11_gen1} is a special case of~\eqref{eq:Pbari_11_gen2}.
Finally, to prove~\ref{item:C12_temp}, consider $z \in \ima Z_{ij} \supseteq \mathcal{S}_i \cap \mathcal{S}_j$. 
By construction, the columns of $Z_{ij}$ are equal to the elements of $\mathcal{R}_{\mathcal{S}_i \cap \mathcal{S}_j}$ (see Section~\ref{sec:preliminaries}). 
Thus,  $z \in \ima Z_{ij}$ if and only if there exist numbers $v_m \in \mathbb{R}$ such that $z = \sum_m v_m r_m$, where the sum is taken over all $r_m \in \mathcal{R}_{\mathcal{S}_i \cap \mathcal{S}_j}$. Now, substitute $z = \sum_m v_m r_m$ into~\eqref{eq:cont_bar}, such that, for all $z\in \ima Z_{ij}$, or equivalently, for all $v_m \in \mathbb{R}$ and all $\hat{x} \in \mathbb{R}^{n-m}$, one finds that
\begin{align}
    x^\top (P_i-P_j) x &= z^\top (\bar{P}_{i}^{11} - \bar{P}_{j}^{11}) z + 2 \hat{x}^\top (\bar{P}_i^{21} - \bar{P}_j^{21}) z \nonumber \\
    &\quad + \hat{x}^\top (\bar{P}_i^{22} - \bar{P}_j^{22}) \hat{x} \nonumber \\
    &= \left(\sum_m v_m r_m\right)^\top (\bar{P}_{i}^{11} - \bar{P}_{j}^{11}) \left(\sum_n v_n r_n\right) \nonumber \\ 
    &\quad + 2 \hat{x}^\top (\bar{P}_i^{21} - \bar{P}_j^{21}) \sum_m v_m r_m \nonumber \\
    &\quad + \hat{x}^\top (\bar{P}_i^{22} - \bar{P}_j^{22}) \hat{x} \nonumber \\
    &= \sum_m \sum_n v_m v_n r_m^\top (\bar{P}_i^{11} - \bar{P}_j^{11}) r_n \nonumber \\
    &\quad + 2\hat{x}^\top \sum_m v_m(\bar{P}_i^{21} - \bar{P}_j^{21}) r_m \nonumber \\
    &\quad + \hat{x}^\top (\bar{P}_i^{22} - \bar{P}_j^{22}) \hat{x} = 0 , \label{eq:Pbari_gen}
\end{align}
where the sums are taken over all $r_m,r_n \in \mathcal{R}_{\mathcal{S}_i \cap \mathcal{S}_j}$, and where the last equality follows from~\eqref{eq:Pbari22_gen},~\eqref{eq:Pbari21_gen}, and~\eqref{eq:Pbari_11_gen2}. Since, by definition, $z=Cx$, \ref{item:C12_temp} follows.

$\mathbf{\ref{item:C12_temp} \iff~\ref{item:C2}}$.
Let $W_{ij}$ be given as in~\ref{item:C2}. 
Consider again the state transformation $\bar{x} = T x$, with $T = \M{C^\top & C_\perp}^\top$, such that $\bar{x} = \M{z^\top & \hat{x}^\top}^\top$, with $z = Cx \in \mathbb{R}^m$ and $\hat{x} = C_\perp^\top x \in \mathbb{R}^{n-m}$. Moreover, $Cx \in \ima Z_{ij}$ if and only if there exists a real vector $v$ such that $Cx = z = Z_{ij}v$. Hence,
\begin{align}
    x &= T^{-1}\bar{x} = \M{C \\ C_\perp^\top}^{-1} \M{z \\ \hat{x}} \nonumber \\
    &=  \M{C \\ C_\perp^\top}^{-1} \M{Z_{ij}v \\ \hat{x}} = \M{C \\ C_\perp^\top}^{-1} \M{Z_{ij} & 0 \\ 0 & I} \M{v \\ \hat{x}} = W_{ij}w ,
\end{align}
for some $w = \M{v^\top & \hat{x}^\top}^\top$ if and only if $Cx \in \ima Z_{ij}$. By substituting $x = W_{ij}w$ into~\eqref{eq:C12_temp}, we obtain equivalently
\begin{equation}
    w^\top W_{ij}^\top (P_i - P_j) W_{ij} w = 0, \textup{ for all real } w,
\end{equation}
which, due to symmetry of $W_{ij}^\top (P_i - P_j) W_{ij}$, is equivalent to~\eqref{eq:C2}.

$\mathbf{\ref{item:C12_temp} \iff~\ref{item:C4}}$. For each pair $(i,j) \in \mathcal{N} \times \mathcal{N}$, consider Corollary~\ref{cor:onesided_projection} with the substitutions
\begin{subequations}
\begin{align}
    \boldsymbol{U} &\coloneq H_{ij} C , \\
    \boldsymbol{Q} &\coloneq P_i - P_j. 
\end{align}
\end{subequations}
With the above substitutions,~\ref{item:OP1} in Corollary~\ref{cor:onesided_projection} reads as follows: There exists a matrix $\boldsymbol{X}\in\mathbb{R}^{m\times n}$ such that 
\begin{equation} \label{eq:5_corr1_1}
    P_i - P_j+(H_{ij} C)^\top \boldsymbol{X} + \boldsymbol{X}^\top (H_{ij} C)= 0,
\end{equation}
i.e., exactly as \ref{item:C4} of Theorem~\ref{thm:1}. On the other hand,~\ref{item:OP2} in Corollary~\ref{cor:onesided_projection} reads
\begin{equation} \label{eq:5_corr1_2}
    x^\top (P_i - P_j)x=0, \text{ for all } x \in \ker H_{ij} C,
\end{equation}
which, because
\begin{equation*}
    x \in \ker H_{ij} C \iff Cx \in \ker H_{ij} = \ima Z_{ij},
\end{equation*}
is equivalent to~\ref{item:C12_temp}. 
Due to Corollary~\ref{cor:onesided_projection},~\eqref{eq:5_corr1_1} is equivalent to~\eqref{eq:5_corr1_2}, and thereby,~$\ref{item:C12_temp} \iff~\ref{item:C4}$.

$\mathbf{\ref{item:C1} \implies~\ref{item:C3}}$.
Suppose that $x^\top (P_i-P_j) x = 0$, for all $Cx \in \mathcal{S}_i \cap \mathcal{S}_j$. As in the proof of $\ref{item:C1} \implies~\ref{item:C12_temp}$, it follows that~\eqref{eq:Pbari22_gen},~\eqref{eq:Pbari21_gen}, and~\eqref{eq:Pbari_11_gen2} hold. Hence, for every $r_m,r_n \in \mathcal{R}_{\mathcal{S}_i \cap \mathcal{S}_j}$
\begin{subequations}\label{eq:phi_elements_new}
\begin{align}
    r_m^\top \bar{P}_{i}^{11} r_n = r_m^\top \bar{P}_{j}^{11} r_n &\eqcolon \phi_{mn} , \label{eq:phi_elements_new_a} \\
    \bar{P}_{i}^{21} r_m = \bar{P}_{j}^{21} r_m &\eqcolon \phi_m , \label{eq:phi_elements_new_b} \\
    \bar{P}_{i}^{22} = \bar{P}_{j}^{22} &\eqcolon \Phi^{22} , \label{eq:phi_elements_new_c}
\end{align}
\end{subequations}
where $\phi_{mn} \in \mathbb{R}$, $\phi_m \in \mathbb{R}^{n-m}$, and $\Phi^{22} \in \mathbb{S}^{n-m}$. Recall that $\mathcal{R}_{\mathcal{S}_i \cap \mathcal{S}_j}$ denotes the set of extremal rays of $\mathcal{S}_i \cap \mathcal{S}_j$.

From~\eqref{eq:phi_elements_new_a}, and on the basis of \cite[Lemma 1]{Ambrosino2015}, there exists a symmetric matrix $\Phi^{11} \coloneq \left\{\phi_{pq}\right\} \in \mathbb{S}^{r}$ for all $p,q \in \{1,2,\ldots, r\}$. Clearly, $\Phi^{11}$ can always be constructed, by collecting the elements in~\eqref{eq:phi_elements_new_a} and giving arbitrary values to the remaining elements (see \cite[Remark 4]{Ambrosino2015}).
On a per region basis, using the extremal ray matrices $R_i$ of ${\mathcal{S}}_i$, one can write
\begin{equation}
R_i ^\top \bar{P}_i^{11} R_i = E_i^\top \Phi^{11} E_i , 
\end{equation}
which follows from the construction of the extraction matrices $E_i$ (see Section~\ref{sec:preliminaries}). Since $R_i$ is invertible, one finds 
\begin{equation}\label{eq:Pbari11_new}
    \bar{P}_i^{11} = R_i^{-\top} E_i^\top \Phi^{11} E_i R_i^{-1} .
\end{equation}
In a similar manner as before, collecting the elements in~\eqref{eq:phi_elements_new_b} in a matrix $\Phi^{21} = \left\{\phi_p\right\} \in \mathbb{R}^{(n-m)\times r}$ for all $p \in \{1,2,\ldots, r\}$, results in
\begin{equation}
    \bar{P}_i^{21} R_i = \Phi^{21} E_i, 
\end{equation}
such that, by invertibility of $R_i$, one finds
\begin{equation}\label{eq:Pbari21_new}
    \bar{P}_i^{21} = \Phi^{21} E_i R_i^{-1} .
\end{equation}
Using~\eqref{eq:phi_elements_new_c},~\eqref{eq:Pbari11_new}, and~\eqref{eq:Pbari21_new}, the partitioned matrix in~\eqref{eq:Pipart_new} is equivalently written as
\begin{equation}\label{eq:Pbari_new}
\begin{split}
    \bar{P}_i& = \M{R_i^{-\top} E_i^\top \Phi^{11} E_i R_i^{-1} & \star \\ \Phi^{21} E_i R_i^{-1} & \Phi^{22}} \\
    &= \M{E_i R_i^{-1} & 0 \\ 0 & I}^\top \M{\Phi^{11} & \star \\ \Phi^{21} & \Phi^{22}}\M{E_i R_i^{-1} & 0 \\ 0 & I}.
\end{split}
\end{equation}

\noindent Then, using $P_i = T^\top \bar{P}_i T$ with $T = \M{C^\top & C_\perp}^\top$, one finds
\begin{align}
    P_i & = \M{C \\ C_\perp^\top}^\top \M{E_i R_i^{-1} & 0 \\ 0 & I}^\top
    \M{\Phi^{11} & \star \\
    \Phi^{21} & \Phi^{22}}
    \M{E_i R_i^{-1} & 0 \\ 0 & I} \M{C \\ C_\perp^\top} \nonumber \\
    & = \M{E_i R_i^{-1} C \\ C_\perp^\top}^\top
    \M{\Phi^{11} & \star\\
    \Phi^{21}  &  \Phi^{22}} \M{E_i R_i^{-1} C \\ C_\perp^\top} . \label{eq:Pi_temp}
\end{align}
By assumption, the matrix $V$ satisfies \magne{$\ima V^\top \supseteq \ker C = \ima C_\perp$}. It follows that there exists a matrix $X$ such that $XV = C_\perp^\top$.
Then, continuing from~\eqref{eq:Pi_temp} with $C_\perp=(XV)^\top$, one finds
\begin{align}
    P_i & = \M{E_i R_i^{-1} C \\ XV}^\top
    \M{\Phi^{11} & \star\\
    \Phi^{21}  &  \Phi^{22}} \M{E_i R_i^{-1} C \\ XV} \nonumber \\
    & = \M{E_i R_i^{-1} C \\ V}^\top
    \M{\Phi^{11} & \star\\
    X^\top\Phi^{21}  &  X^\top\Phi^{22}X} \M{E_i R_i^{-1} C \\ V} \nonumber \\
    &= {F}_i^\top \Phi {F}_i,
\end{align}
where $F_i = \M{(E_i R_i^{-1} C)^\top & V^\top}^\top$.

Hence,~\ref{item:C1} implies that there exists a symmetric matrix $\Phi$ such that $P_i = {F}_i^\top \Phi {F}_i$ for all $i \in \mathcal{N}$.

$\mathbf{\ref{item:C3} \implies~\ref{item:C1}}$.
Let $F_i = \M{(E_i R_i^{-1}C)^\top & V^\top}^\top$, where $V$ satisfies \magne{$\ima V^\top \supseteq \ker C$}. Suppose that there exists a symmetric matrix $\Phi$, such that, $P_i = F_i^\top \Phi F_i$ for all $i\in\mathcal{N}$. Recall from Section~\ref{sec:preliminaries} that the extremal ray matrix of each simplicial cone is constructed as $R_i = \bar{R}E_i$, where $\bar{R} \in \mathbb{R}^{m\times r}$ contains all distinct extremal rays of the simplicial conic partition $\{S_i\}_{i=1}^N$, and $E_i \in \mathbb{R}^{r\times m}$ are selection matrices.
 
Clearly, $Cx \in \mathcal{S}_i \cap \mathcal{S}_j$ if and only if there exist vectors $\lambda_i, \lambda_j \geq 0$ such that $Cx = R_i \lambda_i =  R_j \lambda_j$. 
Furthermore, $E_i \lambda_i = E_j \lambda_j$ if and only if $Cx \in \mathcal{S}_i \cap \mathcal{S}_j$, due to Assumption~\ref{ass:partition}. Hence, for all $Cx \in \mathcal{S}_i \cap \mathcal{S}_j$
\begin{align}
      (E_iR_i^{-1}C - E_jR_j^{-1}C)x &= (E_iR_i^{-1}R_i \lambda_i - E_jR_j^{-1} R_j \lambda_j) \nonumber \\
      &= (E_i \lambda_i - E_j \lambda_j) = 0. \label{eq:Fi_hat_equality}
 \end{align}
As such,
 \begin{align}
     (F_i - F_j)x &= \left( \M{E_i R_i^{-1}C \\ V} - \M{E_j R_j^{-1}C \\ V} \right)x \nonumber \\
     &= \M{(E_i R_i^{-1}C - E_j R_j^{-1}C)x \\ 0} \nonumber \\
     &\stackrel{\eqref{eq:Fi_hat_equality}}{=} 0, \textup{ for all } Cx \in \mathcal{S}_i \cap \mathcal{S}_j. \label{eq:Fi_equality}
 \end{align}
From~\eqref{eq:Fi_equality}, it follows that $F_i x = F_j x$ for all $Cx \in \mathcal{S}_i \cap \mathcal{S}_j$. Hence,
\begin{align}
    x^\top (P_i - P_j) x &= x^\top (F_i^\top \Phi F_i - F_j^\top \Phi F_j)x \nonumber \\
    &= x^\top F_i^\top \Phi F_i x - x^\top F_j^\top \Phi F_j x \nonumber \\
    &\stackrel{\eqref{eq:Fi_equality}}{=} 0, \textup{ for all } Cx \in \mathcal{S}_i \cap \mathcal{S}_j,
\end{align}
i.e.,~\ref{item:C1} is satisfied.

As we have shown $\ref{item:C1} \iff~\ref{item:C12_temp}$, $\ref{item:C12_temp} \iff~\ref{item:C2}$, $\ref{item:C12_temp} \iff~\ref{item:C4}$, and $\ref{item:C1} \iff~\ref{item:C3}$, the proof is complete.$\hfill \square$

\section{Examples}\label{sec:examples}
\begin{example} \label{ex:1} 
Consider the PWQ function $V: \mathcal{Z} \rightarrow \mathbb{R}_{\geq 0}$ given by
\begin{equation} \label{eq:ex1}
    V(x) = \begin{cases}
        x^\top P_1 x , \quad \textup{when } R_1^{-1} x \geq 0, \\
        x^\top P_2 x , \quad \textup{when } R_2^{-1} x \geq 0,
    \end{cases}
\end{equation}
where $\mathcal{Z} = \{ x \in \mathbb{R}^2 \mid x_2 \geq 0 \}$ and
\begin{equation} \label{eq:ex1_matrices}
    P_1 = I = R_1, \; P_2 = \M{2 & -1 \\ -1 & 1}, \; R_2 = \M{-1 & 0 \\ 0 & 1}.
\end{equation}

\begin{figure}
    \centering
    \includegraphics[width=0.9\linewidth]{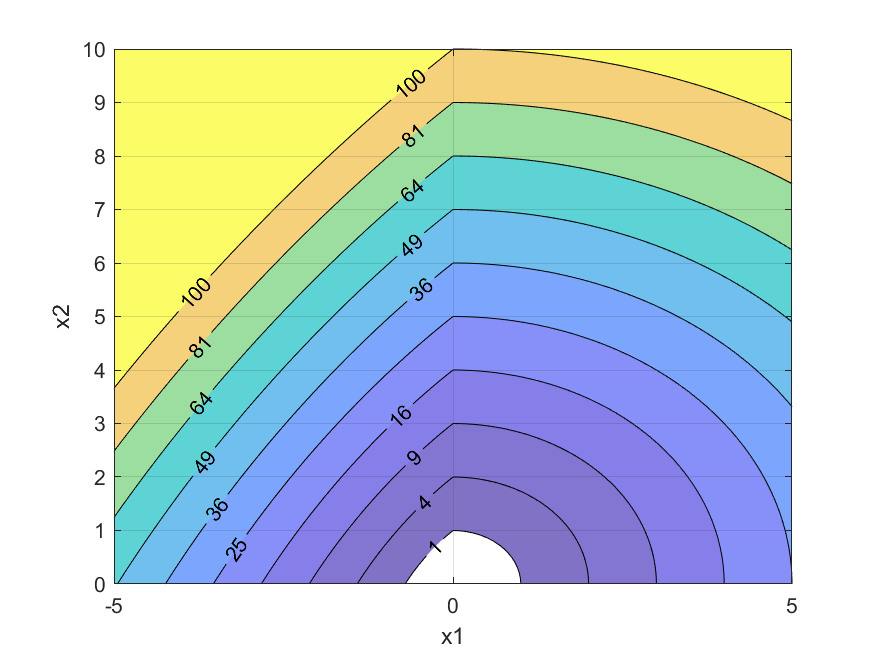}
    \caption{Contour plot of the PWQ function~\eqref{eq:ex1} in example~\ref{ex:1}.}
    \label{fig:ex1}
\end{figure}


\noindent A contour plot of the function $V$ in~\eqref{eq:ex1} is shown in Fig.~\ref{fig:ex1}. Since $C=I$, the matrix $C_\perp$ is omitted. The different continuity conditions in Theorem~\ref{thm:1} will be analysed in detail. 

$\mathbf{\ref{item:C1}}$: The boundary region between $\mathcal{S}_1$ and $\mathcal{S}_2$ is given by $\mathcal{S}_1 \cap \mathcal{S}_2 = \{x\in\mathbb{R}^2 \mid x_1 = 0, x_2\geq 0\}$. At the boundary, we have $V_1(0,x_2) = V_2(0,x_2) = x_2^2$. Hence,~\ref{item:C1} of Theorem~\ref{thm:1} is satisfied and $V(x)$ is (locally Lipschitz) continuous. The continuity of $V(x)$ is also observed in Fig.~\ref{fig:ex1}.

$\mathbf{\ref{item:C12_temp}}$: The so-called image representation of the boundary region is given by $\ima Z_{12} = \{ x\in \mathbb{R}^2 \mid x_1 = 0 \}$ (compared to $\mathcal{S}_1 \cap \mathcal{S}_2$, this set also includes $x_2<0$). Since $(-x)^\top P_1 (-x) = x^\top P_1 x$ holds trivially, it follows that \ref{item:C12_temp} holds.

$\mathbf{\ref{item:C2}}$: With $T=C=I$, everything except the top-left block of $W_{ij}$ in~\eqref{eq:C2}, is omitted. Hence, the condition is reduced to checking whether $Z_{ij}^\top (P_i - P_j) Z_{ij} = 0$ holds. By direct computation, one gets
\begin{equation}
    Z_{12}^\top (P_1 - P_2) Z_{12} = \M{0 \\ 1}^\top \M{1 - 2 & 0 + 1 \\ 0 + 1 & 1 - 1} \M{0 \\ 1} = 0 ,
\end{equation}
where $Z_{12} = \M{0 & 1}^\top$ corresponds to the only boundary region. Hence, $\ref{item:C2}$ is satisfied, as we expected (due to the equivalence).

$\mathbf{\ref{item:C3}}$: Since $C=I$, the matrix $V$ is omitted, such that $P_i = (E_i R_i^{-1})^\top \Phi E_i R_i^{-1}$ and $\Phi \in \mathbb{R}^{3 \times3}$. We order the extremal rays as
\begin{equation}
    \bar{R} = \M{r_1 & r_2 & r_3} = \M{1 & -1 & 0 \\ 0 & 0 & 1} ,
\end{equation}
and define the selection matrices, $E_i$, as
\begin{equation}
    E_1 = \M{1 & 0 \\ 0 & 0 \\ 0 & 1}, \quad E_2 = \M{0 & 0 \\ 1 & 0 \\ 0 & 1}.
\end{equation}
We find that the matrix
\begin{equation}
    \Phi = \M{ 1 & c & 0 \\ c & 2 & 1 \\ 0 & 1 & 1},
\end{equation}
satisfies $P_i = (E_i R_i^{-1})^\top \Phi E_i R_i^{-1}$, for any $c \in \mathbb{R}$, and~\ref{item:C3} holds.

$\mathbf{\ref{item:C4}}$: The matrix $H_{12}$ is constructed as $H_{12} = H_{21} = \M{1 & 0}$. Now, with $\Gamma_{12} \coloneq \M{-1/2 & 1}$, we get
\begin{equation}
    H_{12}^\top \Gamma_{12} + \Gamma_{12}^\top H_{12} = \M{-1 & 1 \\ 1 & 0} = P_1 - P_2,
\end{equation}
and $\Gamma_{21} = -\Gamma_{12}$ follows. Hence,~\ref{item:C4} is satisfied.

\end{example}

\begin{example} \label{ex:ex2}
\magne{
Consider the mass-spring system
\begin{equation}
    m\Ddot{z} + k z = u,
\end{equation}
where $m = 1$ \si{\kilo\gram} is the mass, $k = 1$ \si{\newton\per\meter} is the spring constant, and $u$ is the (force) input. 
Consider the input $u = -f(k_c; z,\dot z) z - f(d_c; z, \dot z) \dot z$ with $k_c = 5$ \si{\newton\per\meter} and $d_c = 2$ \si{\newton\second\per\meter}, where $f(\cdot)$ is a piecewise constant function given by
\begin{equation}
    f(a; z,\dot{z}) = \begin{cases}
        a , \quad \text{when } z\dot{z} \geq 0 \\
        0 , \quad \text{when } z\dot{z} \leq 0.
    \end{cases}
\end{equation}
The idea behind this controller is to apply an opposing force when the mass is moving away from the origin, i.e., when $z\dot z > 0$.
Now, the closed-loop system is 
\begin{equation}
    m\Ddot{z} + f(d_c; z,\dot z) \dot{z} + (k + f(k_c; z,\dot z)) z = 0, 
\end{equation}
i.e., a mass-spring-damper system with variable spring and damping constants.}
\magne{With $x\coloneq \M{z & \dot{z}}^\top$, we obtain the PWL model
\begin{equation} \label{eq:ex2_pwl_system}
    \dot x = A_i x , \quad \textup{when } R_i^{-1} x \geq 0,
\end{equation}
where
\begin{subequations} \label{eq:ex2_system_matrices}
\begin{align}
    A_1 &= A_3 = \M{0 & 1 \\ -\frac{k + k_c}{m} & -\frac{d_c}{m}} = \M{0 & 1 \\ -6 & -2}, \\
    A_2 &= A_4 = \M{0 & 1 \\ -\frac{k}{m} & 0} = \M{0 & 1 \\ -1 & 0}, \\
    R_1 &= - R_3 = \M{1 & 0 \\ 0 & 1}, \quad R_2 = - R_4 = \M{-1 & 0 \\ 0 & 1}.
\end{align}
\end{subequations}}

\begin{figure}
    \centering
    \includegraphics[width=0.9\linewidth]{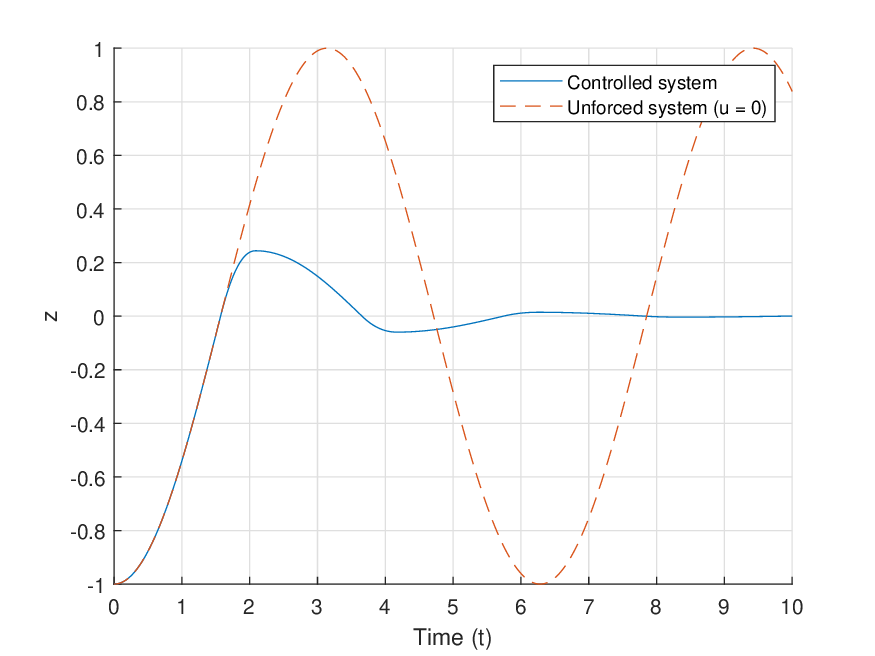}
    \caption{Position $z$ of the system analysed in example~\ref{ex:ex2}, with initial condition $(z(0), \dot z(0)) = (-1, 0)$.}
    \label{fig:ex2}
\end{figure}

\noindent \magne{An example trajectory of the PWL system in \eqref{eq:ex2_pwl_system} is shown in Fig.~\ref{fig:ex2}.
We aim to certify so-called PWQ stability using continuity conditions together with the following LMIs (see, e.g., \cite{Johansson1998, Iervolino15}) 
\begin{subequations} \label{eq:pwq_stability}
\begin{align}
    P_i -R_i^{-\top} W_i R_i^{-1} &\succ 0 , \\
    A_i^\top P_i + P_i A_i + R_i^{-\top} U_i R_i^{-1} &\prec 0 ,
\end{align}
\end{subequations}
where $U_i, W_i \in \mathbb{S}_{\geq 0}^{2}$. Note that~\eqref{eq:pwq_stability} does not guarantee stability of sliding-modes. However, it can be shown that the PWL system analysed in this example does not exhibit any sliding-mode behaviour.}

\magne{The LMIs are solved using MOSEK via the YALMIP MATLAB interface \cite{MOSEK,YALMIP}. We partition each quadrant into two, for a total of eight regions. For ease of representation we restrict ourself to symmetric PWQ functions, so that $P_i = P_{i-4}$ for $i \in \{5,6,7,8\}$. For the condition \ref{item:C3} we enforce this symmetry using that $r_5 = \M{-1 & 0}^\top = - r_1$ and letting $E_4 = \footnotesize \M{0 & 0 & 0 & 1 & 0 & 0 & 0 & 0 \\ -1 & 0 & 0 & 0 & 0 & 0 & 0 & 0}^\top \normalsize$. Then $P_i$, $i\in\{1,\dots,8\}$, with $P_i = P_{i-4}$ for $i \in \{5,6,7,8\}$, defines a continuous PWQ function, and all elements outside the top $4\times 4$ block of $\Phi$ (denoted $\Phi_{11}$) are arbitrary.
With the parametrisation $P_i = F_i^\top \Phi F_i$ in \ref{item:C3}, solving the LMIs in \eqref{eq:pwq_stability} for $i \in \{1,\dots,4\}$ yields
\begin{equation}
    \Phi_{11} =
        \footnotesize
        \def\arraystretch{\arraystretchval}\left[\begin{array}{@{}cc|cc@{}}
            1.5366 & \star & \star & \star \\
            1.1595 & 0.9907 & \star & \star \\\hline
            0 & 0.2935 & 0.3661 & \star \\
            -0.8536 & 0 & 0.4725 & 0.9237
        \end{array}\right],
        \normalsize
\end{equation}
where the entries indicated as zeros are arbitrary. Note that solving \eqref{eq:pwq_stability} only for $i \in \{1,\dots,4\}$ is sufficient, because the PWL system in \eqref{eq:ex2_pwl_system} has symmetric dynamics; that is, the dynamics at any point $x$ are equal to those at $-x$.}

\magne{If, instead, the equality constraints from~\ref{item:C2} are used to ensure continuity, the solver returns (only four digits shown due to space constraints)
\begin{equation} \label{eq:ex2_numerics_equality_constraints}
    \def\arraystretch{\arraystretchval}\left[\begin{array}{@{}c|c@{}}
        P_1 & P_2\\\hline
        P_3 & P_4
    \end{array}\right] =
    \footnotesize
    \def\arraystretch{\arraystretchval}\left[\begin{array}{@{}cc|cc@{}}
        1.6129 & 0.1037 & 1.5905 & 0.0548 \\
        0.1037 & 0.2479 & 0.0548 & 0.3682 \\\hline
        0.9413 & -0.3214 & 1.6129 & 0.3417 \\
        -0.3214 & 0.3682 & 0.3417 & 1.0227
    \end{array}\right].
    \normalsize
\end{equation}
In this case, for example, with $Z_{12} = \frac{1}{\sqrt{2}}\M{1 & 1}^\top$, we obtain $Z_{12}^\top(P_1-P_2)Z_{12} = 5.6\cdot 10^{-12} \neq 0$ when using the full MATLAB-precision matrices, and  $Z_{12}^\top(P_1-P_2)Z_{12} = -~5~\cdot~10^{-5}~\neq~0$ for the four-digit values in~\eqref{eq:ex2_numerics_equality_constraints}, thereby violating the continuity conditions. Note that both the four-digit version in~\eqref{eq:ex2_numerics_equality_constraints} and the full-precision matrices (omitted for brevity) satisfy~\eqref{eq:pwq_stability}.
} 
\end{example}

\section{Conclusion}\label{sec:conclusion}
Continuity conditions known in the literature for PWQ functions on simplicial conic partitions are shown to be equivalent. This result is particularly useful in the context of stability analysis of PWL systems using PWQ Lyapunov functions.
Moreover, the results demonstrate that the method for ensuring continuity solely can be based on practical considerations, without introducing additional conservatism in the analysis.
Finally, a technical lemma useful for showing equivalence of the continuity conditions is presented. This lemma is of independent interest and has significant potential for applications beyond those explored in this technical note.


\bibliographystyle{ieeetr}
\bibliography{references}
\immediate\write18{cp output.aux main.aux}
\end{document}